\documentstyle[12pt,amscd,amssymb]{amsart}
\newtheorem{pr}{Proposition}
\newtheorem{lm}{Lemma}

\newcommand{\proj}{\Bbb P}
\newcommand{\lsys}{\cal{P}}
\newcommand{\sev}{\cal{N}}

\newcommand{\com}{\Bbb C}
\newcommand{\barr}{\overline}

\newcommand{\rarr}{\rightarrow}
\newcommand{\oh}{{\cal{O}}}

\newcommand{\eqq}{\stackrel {\sim}{=}}

\begin{document}
\title{Severi Degrees in Cogenus 3}
\author{J. Harris and R. Pandharipande}
\date{6 April 1995}
\maketitle
\pagestyle{plain}
\setcounter{section}{-1}
\section{{\bf Introduction}}
\subsection{Summary}
Denote by $\lsys (d)$ the linear system of degree $d$ curves in
the complex projective plane $\proj ^2$. $\lsys (d)$ is
a projective space of dimension ${d+2 \choose 2}-1$. Let
$\sev(n,d) \subset \lsys (d)$ be the subset corresponding to
reduced, nodal curves with exactly $n$ nodes. $\sev(n,d)$ is
empty unless $0\leq n \leq {d \choose 2}$.
Points of
$\sev(n,d)$ may correspond to reducible curves.
If nonempty,
$\sev(n,d)$ is a quasi-projective subvariety of
pure codimension $n$ in $\lsys (d)$. Let $\barr{\sev}(n,d)$ denote
the closure of $\sev(n,d)$ in $\lsys (d)$. In this paper,
formulas for the
the degree of $\barr{\sev}(n,d)$ for $n=1,2,3$ are computed:
\begin{equation}
\label{one}
f_1(d)= 3(d-1)^2,
\end{equation}
\begin{equation}
\label{two}
f_2(d)= {3 \over 2}  (d-1)(d-2)(3d^2-3d-11),
\end{equation}
\begin{equation}
\label{three}
f_3(d)= {9 \over 2} d^6 - 27d^5+
          {9 \over 2} d^4 +
          {423 \over 2 } d^3 - 229 d^2 -
         {829 \over 2 } d +525,
\end{equation}

$$\forall d \geq 1, \ \ degree(\barr{\sev}(1,d))=f_1(d)$$
$$\forall d\geq 3, \ \ degree(\barr{\sev}(2,d))= f_2(d)$$
$$\forall d\geq 3, \ \ degree(\barr{\sev}(3,d))= f_3(d).$$

These formulas are classical. The computation presented here
is new.
The method involves the geometry of the Hilbert scheme
of points in $\proj^2$ and the Bott residue formula (following the
technique developed in [E-S]).
The most successful methods for obtaining cogenus formulas appear
in [V]. Via a sophisticated singularity analysis, I. Vainsencher
obtains the above results and the following further cogenus formulas:
\begin{eqnarray*}
f_4(d) & = & {27\over 8}d^8-27 d^7+{1809\over 4} d^5-642 d^4-2529 d^3
+{37881\over 8} d^2+{18057\over 4}d-8865, \\
f_5(d) & = & {81\over 40} d^{10}- {81\over 4} d^9 - {27\over 8}d^8+
{2349\over 4}d^7-1044 d^6-{127071\over 20} d^5 +{128859\over 8} d^4  \\
& & +{59097\over 2}d^3-{3528381\over 40}d^2 -{946929\over 20}d+153513, \\
f_6(d) & = & {81\over 80} d^{12}-{243\over 20}d^{11} -{81\over 20} d^{10}
+{8667\over 16}d^9 -{9297\over 8} d^8 -{47727\over 5} d^7 \\
& & + {2458629\over 80} d^6
 +{3243249\over 40}d^5 -{6577679\over 20}d^4-{25387481\over 80} d^3
+{6352577\over 4} d^2 \\ & & +{8290623\over 20}d -2699706. \\
\end{eqnarray*}

The idea for the computation presented here originated in a
conversation at the October 1994 Utah-UCLA-Chicago Algebraic
Geometry Workshop.
Discussions with W. Fulton were helpful.
The second author benefitted from conversations with
D. Edidin and W. Graham. Thanks are due to S. Stromme for pointing
out I. Vainsencher's results.

\subsection{The Method}
\label{method}
For each $n\geq 1$,
let $H(n)$ be the Hilbert scheme of length $n$ subschemes of $\proj ^2$.
$H(n)$ is a nonsingular variety of dimension $2n$ with generic element
corresponding to a subscheme of $n$ distinct points of $\proj ^2$.
There exists a rational map
$$\psi_n: H(n) \ - \ - \rarr H(3n)$$
given by squaring the ideal sheaf. For $1\leq n \leq 3$,
we will consider resolutions
$$\barr{\psi}_n: X(n) \rarr H(3n)$$ of $\psi_n$.
It is easily checked that $\psi_1$ and $\psi_2$ are everywhere defined.
Let $X(1)=H(1)$, $\barr{\psi}_1=\psi_1$ and $X(2)=H(2)$,
$\barr{\psi}_2=\psi_2$.
Let $F \hookrightarrow H(3)$ be the locus of
length $3$ subschemes isomorphic to $\com[x,y]/(x^2, xy, y^2)$.
$F$ is a nonsingular subvariety abstractly isomorphic to $\proj^2$.
Since $(x^2,xy,y^2)^2$ is an ideal of length $10$, $\psi_3$ is not
defined on $F$. Let
$X(3)$ be the blow-up of $H(3)$ along $F$. In section
(\ref{resmap}), it is shown that
$\psi_3$ is defined on $H(3) \setminus F$ and extends to a morphism
$$\barr{\psi}_3: X(3) \rarr H(9).$$ The degrees of $\barr{\sev}(1,d)$,
$\barr{\sev}(2,d)$, and $\barr{\sev}(3,d)$ will be expressed as Chern
classes of certain tautological bundles over $X(1)$, $X(2)$, and $X(3)$
respectively.

Let $U(n)\hookrightarrow H(n) \times \proj^2 $ be the universal subscheme
over $H(n)$. For $1\leq n \leq 3$, let
$$Y(n)= X(n) \times_{H(3n)} U(3n) \hookrightarrow X(n) \times \proj^2.$$
Let $\pi_n$, $\rho_n$ be the projections from $Y(n)$ to $X(n)$, $\proj ^2$
respectively. For pairs $(n,d)$ where  $1\leq n \leq 3$ and $\sev(n,d)\neq
\emptyset$,
let
$$E(n,d)= \pi_{n*} \rho_n^*(\oh_{\proj^2}(d)).$$
$E(n,d)$ is easily seen to be a rank $3n$ vector bundle on $X(n)$.
We claim $$c_{2n}(E(n,d))= degree(\barr{\sev}(n,d)).$$
A sketch of the argument is as follows. Let
$1\leq n \leq 3$ and let $d$ be such that $\sev(n,d)\neq \emptyset$.
Since $\barr{\sev}(n,d)$ is of codimension
$n$ in $\lsys (d)$,
the degree is the cardinality of a generic $n$-plane slice. An $n$-plane,
$\cal{L}\subset \lsys (d)$,
is equivalent to an $n+1$-dimensional linear subspace of
$L\subset H^0(\proj^2,\oh_{\proj^2}(d))$. $L$ canonically yields
an $n+1$ dimensional subspace
$\barr{L}\subset H^0(X(n),E(n,d))$.
A generic point $\xi \in X(n)$ corresponds to $n$ points of
$\proj^2$. It is checked that $\barr{L}$ drops rank at
$\xi$ if and only if there exists an element of $\cal{L}$
corresponding to a plane curve singular at the
$n$ points of $\xi$. It is further checked, for generic $\cal{L}$, the singular
plane curve must be reduced, nodal with exactly $n$ nodes. The
nongeneric points of $X(n)$ make no contribution to
the degeneracy locus. The degree of
$\barr{\sev}(n,d)$ is thus equal to the cardinality  of the degeneracy locus of
$n+1$ sections of $E(n,d)$. The latter is the $c_{2n}(E(n,d))$.
Section (\ref{loci}) contains the full argument. In principle,
this approach may be attempted for $n\geq4$. Explicit resolutions
of $\psi_n$ are needed. For $n\geq 4$, a correction term to $c_{2n}(E(n,d))$
for small $d$ is required to account for the contribution of nonreduced
curves to the degeneracy locus.

It remains to compute $c_{2n}(E(n,d))$. There is a diagonal torus action
on $\proj^2$ which can be lifted to $X(n)$ and $E(n,d)$. The Bott
residue formula expresses the desired Chern class in terms of the differential
data of the torus action at fixed points. This approach yields the
degree formulas. It should be mentioned there are more direct ways
of obtaining formulas (\ref{one}) and (\ref{two}).
Since the residue calculations
for formula (\ref{three}) contain those required for
(\ref{one}) and (\ref{two}), we present
a unified approach. The explicit residue computations are presented
in section (\ref{bott})

\section{\bf $\psi_n$ And $\barr{\psi}_n$ For $n=1,2,3$}
\label{resmap}
Consider the universal subschemes $U(n) \hookrightarrow H(n) \times \proj^2$.
Let $U^2(n)$ be subscheme of $H(n)\times \proj^2$ defined by the
square of the ideal of $U(n)$. Let $G(n)\subset H(n)$ denote the generic
locus corresponding to subschemes of  $n$ distinct points of $\proj^2$.
$U^2(n)$ is flat over $G(n)$ of degree $3n$. Therefore there is map
$$\psi_n: G(n) \rarr H(3n).$$
Certainly $G(1)=H(1)$. For $n=2$, the complement of $G(2)$ in $H(2)$ consists
of linear double points (with ideals isomorphic to $(x,y^2)\subset \com[x,y]$).
Since $(x,y^2)^2=(x^2, xy^2,y^4)$ has length $6$, $U^2(2)$ is flat over
$H(2)$. The map $\psi_2$ extends to $H(2)$.

For $n=3$, the situation is more complex. By the results for $n=2$,
$\psi_3$ extends to all of $H(3)$ with the possible exception of the
the triple points. The isomorphism classes of length $3$ subschemes of
$\com^2$ supported at a point are given by the following ideals:
\begin{enumerate}
\item[(i.)] $(x,y^3)\subset \com[x,y]$.
\item[(ii.)] $(x+y^2, x^2, xy) \subset \com[x,y]$.
\item[(iii.)] $(x^2,xy,y^2) \subset \com[x,y]$.
\end{enumerate}
It is easy to check the squares of the ideas of type (i) and (ii) have
length 9. $(x^2,xy,y^2)^2=(x^4,x^3y,x^2y^2,xy^3,y^4)$ has length 10.
Let $F\hookrightarrow H(3)$ be the nonsingular subscheme corresponding to the
points of type (iii). $U^2(3)$ is flat over $H(3) \setminus F$. Therefore
$\psi_3$ extends to
$\psi_3: H(3)\setminus F \rarr H(9).$

Let $V\subset \proj^2$ be a coordinate affine chart.
$V$ is a two dimensional complex vector space.
Let
$$A= \bigoplus_{k=0}^{\infty} Sym^k(V^*)$$
be the affine coordinate ring of $V$.
Let $m\subset A$ be the maximal ideal corresponding to the point
$0\in V$. The ideal $m^2$ is  of type (iii). Let
$[A/m^2]\in H(3)$ denote the Hilbert point corresponding to
$A/m^2$.
The tangent space to $H(3)$ at $[A/m^2]$ is
canonically isomorphic the module of $A$-homomorphisms
$Hom_A(m^2/m^4, A/m^2)$. It is easily seen there are canonical
isomorphims
$$Hom_A(m^2/m^4, A/m^2) \eqq Hom_{\com}(m^2/m^3, m/m^2) \eqq
Hom_{\com}(Sym^2(V^*), V^*).$$
$V$ is canonically identified with the space of invariant
vector fields on $V$. Therefore, there is a canonical map
$$\mu: V \rarr Hom_{\com}(Sym^2(V^*), V^*)$$
given by differentiation of functions. Certainly, $[A/m^2]\in F$.
The tangent space to $F$ at $[A/m^2]$ is canonically isomorphic
to $V$. The map $\mu$ is
the differential
of the inclusion of $F$ in $H(3)$ at $[A/m^2]$.
Consider the exact sequence
\begin{equation*}
0\rarr K \rarr Sym^2(Sym^2(V^*)) \rarr Sym^4(V^*) \rarr 0
\end{equation*}
given by multiplication. $K$ is a one dimension $\com$-vector space.
There is a canonical map $\nu$ and an exact sequence:
\begin{equation}
\label{xact}
0 \rarr V \stackrel{\mu}{\rarr} Hom_{\com}(Sym^2(V^*), V^*)
\stackrel{\nu}{\rarr}  Hom_{\com}(K, Sym^3(V^*)) \rarr 0.
\end{equation}
Briefly, an element of $\gamma \in Hom_{\com}(Sym^2(V^*), V^*)$
yields a map $id+ \gamma: Sym^2(V^*) \rarr A$.
Multiplication of $id+\gamma$ induces
a map
$Sym^2(Sym^2(V^*)) \rarr A$.
The latter map takes $K$ to $Sym^3(V^*)$. The exactness of (\ref{xact})
is a simple exercise.

Let $X(3)$ be the blow up of $H(3)$ along $F$.
By sequence (\ref{xact}), there is a natural correspondence between the fiber
of the
projective normal bundle of  $F$ in $H(3)$ at $[A/m^2]$ and the
projective space $\proj( Hom_{\com}(K, Sym^3(V^*)))$. The map
$\psi_3$ can be extended to the projective normal bundle of $F$
by mapping an element $[\xi] \in \proj(Hom_{\com}(K,Sym^3(V^*)))$ to
the ideal of length $9$ given by $(m^4, image(\xi))\subset A$.
We have defined a map
$$\barr{\psi}_3: X(3) \rarr H(9).$$
It is not hard to check that $\barr{\psi}_3$ is an algebraic morphism.

\section{\bf $E(n,d)$, Degeneracy Loci, and $degree(\barr{\sev}(n,d))$.}
\label{loci}
Let $n=1,2,$ or $3$. Let $d$ be an integer such that
$\sev(n,d)$ is nonempty. Following the notation of section (\ref{method}),
 $$E(n,d)=\pi_{n*} \rho_n^*(\oh_{\proj^2}(d)).$$
Let $\cal{L}\subset \lsys (d)$ be an $n$-plane corresponding
to an $n+1$ dimensional subspace
$L\subset H^0(\proj^2,\oh_{\proj^2}(d))$.
Let $\barr{L}$ denote the naturally induced $n+1$ dimensional subspace
of $H^0(X(n), E(n,d))$. For any $\xi\in X(n)$, let
$\barr{L}_{\xi} \subset E(n,d)_{\xi}$ be the subspace of the fiber
generated by $L$. Let $I_{\xi}\subset \oh_{\proj^2}$ denote the ideal sheaf of
subscheme corresponding to $\barr{\psi}_n(\xi)$.

\begin{lm}
\label{taut}
Let $[\xi]\in X(n)$. Then, $dim(\barr{L}_{\xi}) < n+1$ if and only if
there exists a nonzero $l \in L$ such that $l\in H^0(\proj^2, I_{\xi}(d))$.
\end{lm}
\noindent Lemma (\ref{taut}) is a tautological statement. We abuse notation
slightly to let $G(n)\subset X(n)$ denote the locus of points $\xi \in X(n)$
such that $\barr{\psi}_n(\xi)$ corresponds to a subscheme of $n$ distinct
points
of $\proj^2$.

\begin{lm}
\label{nobound}
For generic $\cal{L}$, the following holds:
$$\forall \xi \in X(n)\setminus G(n), \ \ dim(\barr{L}_{\xi})=n+1.$$
\end{lm}
\begin{pf}
For $n=1$, the result is vacuous.

For $n=2$, $d$ must be at least $3$ to ensure $\sev(2,d) \neq \emptyset$.
$X(2)\setminus G(2)$ consists of the $3$ dimensional locus of
linear double points (with ideals isomorphic to $(x, y^2)\subset \com[x,y]$).
The squared ideal $(x^2, xy^2, y^4)$ imposes 6 conditions on
linear series of degree $d \geq 3$. Hence,  the locus of elements
$[l] \in \lsys (d)$ such that $l \in H^0(\proj^2, I_{\xi}(d))$ for some
$\xi \in X(2) \setminus G(2)$ is of codimension at least $6-3=3$. A generic
$2$-plane
has empty intersection with this locus.

For $n=3$, $d$ must again be at least $3$ to ensure $\sev(3,d)\neq \emptyset$.
$X(3)\setminus G(3)$ consists of three loci:
\begin{enumerate}
\item [(a.)] The projective normal bundle of $F$, $\proj(N_F)$.
\item [(b.)] The quasi-projective locus, $B$,  corresponding to ideals of type
(i) and (ii) of section (\ref{resmap}).
\item [(c.)] The quasi-projective locus, $C$, corresponding to
a point $p$  union a linear double point supported at $q\neq p$.
\end{enumerate}
The ideal $I_{\xi}$ for any $\xi$ in $\proj(N_F)$
certainly imposes $9$ conditions on linear series of degree $d\geq 3$.
The dimensions of $\proj(N_F)$  is $5$.
The locus of elements $[l] \in \lsys (d)$
such that $l\in H^0(\proj^2, I_{\xi}(d))$
for some $\xi \in \proj(N_F)$ is codimension at least $9-5=4>3$.

$B$ consists of a $3$ dimensional locus of ideals of type (i) and
a $4$ dimensional locus of ideals of type (ii). Ideals of
type (i) and (ii) are easily seen to impose $7$ and $9$ conditions
respectively
on linear series of degree $d\geq 3$. We see $7-3=4>3$ and
$9-4=5>3$.

If $J$ is the ideal of a linear double point supported at $q\in \proj^2$,
$J$ imposes $6$ conditions
on linear series of degree $d\geq 3$. The generic element of
$H^0(\proj^2, J(d))$ is singular only at $q$. The condition that
$l\in H^0(\proj^2, J(d))$ be singular at some point in $\proj^2 \setminus \{ q
\}$ is
therefore a codimension $1$ condition on $H^0(\proj^2, J(d))$.
The locus of elements $[l] \in \lsys (d)$ such that $l\in H^0(\proj^2,
I_{\xi}(d))$ for some $\xi \in T$ is therefore of codimension at least
$7-3=4>3$. A generic $3$-plane in $\lsys (d)$ avoids the loci
corresponding to subsets (a), (b), and (c) of $X(3) \setminus G(3)$.
\end{pf}

\begin{lm}
\label{nsing}
Let $Sing(n,d)\subset \lsys(d)$ be the quasi-projective locus of
reduced curves with at least $n$ singular points. Then
$Sing(n,d)\subset \barr{\sev}(n,d)$.

\noindent
(This result holds for all $n$.)
\end{lm}
\begin{pf}
Let $[C]\in Sing(n,d)$. We must show $[C]\in \barr{\sev}(n,d)$.
Since $\barr{\sev}(m,d)\subset \barr{\sev}(n,d)$ for $m>n$,
we can reduce to the case were $C$
has exactly $n$ nodes. Let $\sum$ be the $n$ singular points of $C$.
The projective tangent space to $Sing(n,d)$ at
$[C]$ is given by the linear system of
degree $d$ curves passing through $\sum$.
By a study of the adjoint conditions ([ACGH], p.60), $\sum$ imposes
$n$ independent conditions on degree $d$ curves. It also follows from
the adjoint analysis that any subideal of $I(\sum)$ of index $1$ imposes
independent conditions on degree $d$ curves.
Since the condition
of having exactly $n$ nodes is open in $Sing(n,d)$, $Sing(n,d)$ is
nonsingular of codimension $n$ at $[C]$.

        Let $f(x,y)$ be the equation of a plane curve singularity at the
origin $(0,0)$. Consider first order, {\em equisingular}
deformations of the type
$f(x,y)+\epsilon \cdot g(x,y)$. It is a fact ([fact]) that such $g(x,y)$
generate the the maximal ideal of $(0,0)$ if and only if
the singularity of $f(x,y)$ at the origin is a node.

        Suppose $[C]$ is not nodal. The equisingular deformations
of $[C]$ correspond at most to the linear system of degree $d$
curves passing through a subideal of $I(\sum)$ of index $1$. Hence
the these equisingular deformations are of codimension at least $n+1$.
Therefore, the generic member of each component must be nodal.
\end{pf}

\noindent
Finally, a simple dimension analysis yields:
\begin{lm}
\label{reduced}
Let $Nonred(d)\subset \lsys (d)$ be the closed locus of nonreduced curves.
For all $d$, the codimension of $Nonred(d)$ is greater than $2$.
For $d\geq 3$, the codimension of $Nonred(d)$ is greater than $3$.
\end{lm}

Let $\cal{L}\subset \lsys(d)$ be a generic $n$-plane.
By Lemma (\ref{nsing}) the following sequence of inclusions hold:
$$\sev(n,d) \subset Sing(n,d) \subset \barr{\sev}(n,d).$$
Therefore, we obtain
\begin{equation}
\label{nodes}
\cal{L}\cap \sev(n,d)= \cal{L}\cap Sing(n,d)=
\cal{L}\cap \barr{\sev}(n,d).
\end{equation}
The intersection is transverse and,
\begin{equation*}
|\cal{L}\cap {\sev}(n,d) | = degree(\barr{\sev}(n,d)).
\end{equation*}
The degeneracy locus $D$ of $\barr{L}$ is the subscheme where
$\barr{L}$ drops rank. By Lemma (\ref{nobound}), $D$ is supported
in $G(n)$. If $n=1,2,$ or $3$ and $d$ is such that $\sev(n,d)\neq \emptyset$,
$E(n,d)$ is generated on $G(n)$
by sections induced from
$H^0(\proj^2, \oh_{\proj^2}(d))$. (Note $E(n,d)$ need not be generated
by these sections on $X(n)$.) Therefore, the canonical sections
yield a map $$\lambda_n(d): G(n) \rarr Grassmanian(n,d).$$
The degeneracy locus $D$ on $G(n)$ is expressed
as the $\lambda_n(d)$-intersection with a Schubert class determined
by $\barr{L}$. By Kleiman's result on
intersections in homogeneous spaces ([H], p.273),
$D$ is a reduced dimension zero
subscheme of $G(n)$
for generic $\cal{L}$.
If $\xi \in D\subset G(n)$ is a point, by Lemma (\ref{taut}) there is
a curve $[l]\in \cal{L}$ singular at the $n$ points of $\proj^2$
corresponding to $\xi$.
$[l]$ must be reduced by Lemma (\ref{reduced}).
$[l]$ therefore corresponds to
a point of the intersection (\ref{nodes}) and must be nodal with exactly
$n$ nodes.
$[l]$ must be unique or else there would be a pencil of
curves of $\cal{L}$ singular at the $n$ points corresponding to $\xi$.
We have defined a set map $\tau: D \rarr \cal{L}\cap \sev(n,d)$. Since
$D$ can be recovered from the nodes of $[l]$, $\tau$ is injective.
By Lemma (\ref{taut}), $\tau$ is surjective.
Hence
$ |D| =| \cal{L} \cap \sev(n,d)|$.
By the Thom-Porteous formula, $|D|= c_{2n}(E(n,d))$.
\begin{pr}
\label{chern}
For $n=1,2,3$,
$degree(\barr{\sev}(n,d))= c_{2n}(E(n,d))$.
\end{pr}

\section{\bf The Bott Residue Computation}
\label{bott}
\subsection{The Formula}
We first state the form of the Bott Residue Formula ([B]) that will be used.
Let $M$ be a nonsingular variety of dimension $m$
with an algebraic $\com^*$-action.
Let $q\in M$ be a fixed point of the $\com^*$-action.
The differential of the action naturally induces a $\com^*$-representation
on the tangent space $T_q(M)$. Let $\alpha_1(q), \ldots, \alpha_m(q)$
be the $m$
weights of the $\com^*$-representation on $T_q(M)$.
Let $\cal{Q}\subset M$ be the fixed point set. Assume
\begin{enumerate}
\item[(1.)] $\cal{Q}$ is discrete.
\item[(2.)] $\forall q \in \cal{Q}$ and $\forall j$, $\alpha_j(q) \neq 0$.
\end{enumerate}
In fact, condition (2) is a consequence of condition (1).
Suppose $E$ is an algebraic vector bundle of rank $r$ on $M$  with an
equivarient $\com^*$-action. For each $q\in \cal {Q}$, there
is a $\com^*$-representation on $E_q$. Let
$\beta_1(q)$, $\ldots$, $\beta_r(q)$ be the weights
of this $\com^*$-representation.
Finally, let $\sigma_{i,j}(x_1, x_2, \ldots, x_j)$ be the $i^{th}$ elementary
symmetric polynomial in the variables $x_1$, $x_2$, $\ldots$, $x_j$.
If $i>j$, $\sigma_{i,j}=0$.
The Bott Residue Formula expresses $c_m(E)$ in terms of the
$\com^*$-weights at the fixed points:
$$c_m(E)= \sum_{q\in \cal{Q}} {\sigma_{m,r}(\beta_1(q),
\ldots, \beta_r(q)) \over \sigma_{m,m}(\alpha_1(q),
\ldots, \alpha_m(q))}.$$
Since $\sigma_{m,m}$ is the product monomial and $\alpha_j(q)\neq 0$,
the right hand side is well defined.

\subsection{Torus Actions}
Let $Z$ be a $3$ dimensional $\com$-vector with basis
$\barr{z}=(z_0,z_1,z_2)$. Let $\barr{w}=(w_0,w_1,w_2)$ be a
triple  of integral weights. Let
$\lambda(\barr{w})$ be the $\com^*$-representation with
weights $(w_0,w_1,w_2)$ diagonal with respect to $\barr{z}$.
Let $\proj^2= \proj(Z)$.
The representation $\lambda(\barr{w})$ induces a $\com^*$-action
on $\proj(Z)$. There is an induced $\com^*$-action on $H(n)$,
$X(n)$, and $Y(n)$. Recall the natural isomorphism
$H^0(\proj(Z), \oh_{\proj(Z)}(d)) \eqq Sym^d(Z^*)$.
There is a canonical equivariant lifting of $\lambda(\barr{w})$ to
$\oh_{\proj(Z)}(d)$ such that the induced representation on global
sections is $Sym^d(\lambda^*(\barr{w}))$.
The canonical equivariant lifting of $\lambda(\barr{w})$ to
$\oh_{\proj(Z)}(d)$ induces an equivariant $\com^*$-action
on $E(n,d)= \pi_{n*}\rho^*(\oh_{\proj(Z)}(d))$ over $X(n)$.

\subsection{The Case $n=1$}
$X(1)=H(1)=\proj(Z)$.
There are $3$ fixed points for distinct weights $\barr{w}=(w_0,w_1,w_2)$.
Analysis of the fixed point $[z_0]$ yields
\begin{equation*}
\begin{array}{ll}
\alpha_1([z_0])=& w_1-w_0 \\
\alpha_2([z_0])= &w_2-w_0
\end{array}
\end{equation*}
To simplify notation, let $Z_0, Z_1, Z_2$ be a basis of
$Z^*$ dual to $\barr{z}$.
To calculate the action on $E(1,d)$, observe that
$\barr{\psi}_1([z_0])$ is the subscheme given by the ideal
$(Z_1^2, Z_1Z_2, Z_2^2)$. Therefore the sections
$Z_0^d$, $Z_0^{d-1}Z_1$, $Z_0^{d-1}Z_2$ generate the
fibre of $E(1,d)$ at $[z_0]$. We obtain:
\begin{equation*}
\begin{array}{ll}
\beta_1([z_0])= &-dw_0 \\
\beta_2([z_0])= &-(d-1)w_0-w_1 \\
\beta_3([z_0])= & -(d-1)w_0-w_2
\end{array}
\end{equation*}
The analysis for $[z_1]$ and $[z_2]$  is similar.
The Bott Residue Formula now yields $c_2(E(1,d))= 3(d-1)^2$.
Since the Chern class does not depend upon the weights, it
is simplest to fix values $\barr{w}=(0,1,2)$ when using the
residue formula.

\subsection{The Case $n=2$}
\label{case2}
$X(2)=H(2)$. There are $9$ fixed points for distinct weights
$\barr{w}=(w_0,w_1,w_2)$. Three fixed points correspond to
the subschemes:
$$[z_0]\cup [z_1], \ \ [z_0] \cup [z_2], \ \ [z_1]\cup [z_2].$$
For these points, the analysis of the $n=1$ suffices to yield
the $\alpha$ and $\beta$-weights.
There are $6$ fixed points given by the subschemes $D_{ij}=
\oh_{\proj(Z)}/ (Z_i^2, Z_j)$ for ordered pairs $1 \leq
i \neq j \leq 3$.
We carry out the analysis
at the fixed point $[D_{1,2}]$.
In the affine open $Z_0\neq 0$, let $I=((Z_1/Z_0)^2, (Z_2/Z_0))$.
Let $A=\com[(Z_1/Z_0), (Z_2/Z_0)]$.
The tangent
space to $H(2)$ at $[D_{1,2}]$ is canonically isomorphic to
$Hom_A(I/I^2,A/I)$. $I/I^2$ is the free $A/I$ module with
generator $(Z_1/Z_0)^2$ and $(Z_2/Z_0)$. $A/I$ is generated
by $1$ and $(Z_1/Z_0)$. Since we know the $\com^*$-action on
basis elements of $I/I^2$ and $A/I$, we obtain:
\begin{equation*}
\begin{array}{ll}
\alpha_1([D_{1,2}])= & 2w_1-2w_0 \\
\alpha_2([D_{1,2}])=& w_1-w_0 \\
\alpha_3([D_{1,2}])=& w_2-w_0 \\
\alpha_4([D_{1,2}])=& w_2-w_1
\end{array}
\end{equation*}
Note
$\barr{\psi}_2([D_{1,2}])$ is the
subscheme defined by $(Z_1^4, Z_1^2Z_2,Z_2^2)$.
Therefore the elements $Z_0^d$, $Z_0^{d-1}Z_1$, $Z_0^{d-1}Z_2$,
$Z_0^{d-2}Z_1^2$, $Z_0^{d-2}Z_1Z_2$, and $Z_0^{d-3}Z_1^3$ yield
a basis the fiber of $E(2,d)$ over $[D_{1,2}]$.
The $\beta$-weights are therefore:
\begin{equation*}
\begin{array}{ll}
\beta_1([D_{1,2}])= & -dw_0 \\
\beta_2([D_{1,2}])= & -(d-1)w_0-w_1 \\
\beta_3([D_{1,2}])= &-(d-1)w_0-w_2 \\
\beta_4([D_{1,2}])= &-(d-2)w_0-2w_1 \\
\beta_5([D_{1,2}])= &-(d-2)w_0-w_1-w_2 \\
\beta_6([D_{1,2}])= &-(d-3)-3w_1
\end{array}
\end{equation*}
The $\alpha$ and $\beta$-weights at the other points $[D_{i,j}]$
are obtained by appropriate permutations of $w_0$, $w_1$ and
$w_2$ in the above formulas.
After some algebra (MAPLE was used at this point), the Bott
Residue formula yields
$$c_4(E(2,d))={3\over 2} (d-1)(d-2)(3d^2-3d-11).$$

\subsection{The Case n=3}
Consider the $\com^*$-action for weights $\barr{w}$ on
$X(3)$. The analysis at nontriple points reduces to previous
computations.
There is one fixed point corresponding to the subscheme $[z_0]\cup[z_1]
\cup [z_2]$. There are 12 fixed point of the type
$$[z_i]\cup D_{j,k}$$
where $Supp(D_{j,k})\neq [z_i]$. The $\alpha$ and $\beta$-weights
for these 1+12 points are easily  obtained from the weights in the
$n=1$ and $2$ cases.

Next consider fixed points points $\xi \in X(3)\setminus \proj_F(N)$
where $\barr{\psi}_3(\xi)$ is
a triple point.
If $(w_0,w_1,w_2)$ are distinct and no two weights sum to
twice the third, then
there are $6$ such $\xi$. They are given by the subschemes $T_{i,j}=
\oh_{\proj(Z)}/ (Z_i^3,Z_j)$ for $1\leq i \neq j \leq 3$.
We carry out the weight analysis at the point $[T_{1,2}]$.
Following the notation of section (\ref{case2}), let
$A=\com[(Z_1/Z_0), (Z_2/Z_0)]$ be the affine coordinate ring
for $Z_0\neq 0$. Let $I=((Z_1/Z_0)^3, (Z_2/Z_0))$. Since
$X(3)$ is isomorphic to $H(3)$ at $[T_{1,2}]$, the tangent space
is given by $Hom_A(I/I^2, A/I)$. $I/I^2$ is seen to be a free
$A/I$ module of rank $2$ with basis $(Z_1/Z_0)^3$ and $(Z_2/Z_0)$.
We hence obtain the six $\alpha$-weights:
$$3w_1-3w_0,\ \ \ 2w_1-2w_0,\ \ \ w_1-w_0,\ \ \ w_2-w_0,\ \ \ w_2-w_1,\
\ \ w_2+w_0-2w_1.$$
Since $(Z_1^3,Z_2)^2=(Z_1^6, Z_1^3Z_2, Z_2^2)$,
the fiber of $E(3,d)$ at $[T_{1,2}]$
is spanned by the (possibly rational) sections
$Z_0^d$, $Z_0^{d-1}Z_1$, $Z_0^{d-1}Z_2$, $Z_0^{d-2}Z_1^2$,
$Z_0^{d-2}Z_1Z_2$, $Z_0^{d-3}Z_1^3$, $Z_0^{d-3}Z_1^2Z_2$,
$Z_0^{d-4}Z_1^4$, $Z_0^{d-5}Z_1^5$.  Therefore, the nine  $\beta$-weights
are:
\begin{equation*}
\begin{array}{lll}
-dw_0 & -(d-1)w_0-w_1 &  -(d-1)w_0-w_2 \\
-(d-2)w_0-2w_1 & -(d-2)w_0-w_1-w_2 & -(d-3)w_0-3w_1 \\
-(d-3)w_0-2w_1-w_2 & -(d-4)w_0-4w_1 & -(d-5)w_0-5w_1
\end{array}
\end{equation*}
Again, the weights at the other $[T_{i,j}]$ are obtained
by appropriate permutations of $\barr{w}$.

Finally, consider the fixed points $\xi \in \proj_F(N)\subset X(3)$.
Let $m_i$ be the ideal of the point $[z_i]$.
Let $F_i$ for $1\leq i \leq 3$ be the
subscheme $\oh_{\proj(Z)}/ m_i^2$.
$[F_i]\in F$. The points $[F_i]$ are the 3 fixed points of the
$\com^*$-action on $F$. The fixed points of $\proj_F(N)$ must
lie in the fibers of $\proj_F(N)$ over the $[F_i]$.
We analyze the case $i=1$. The fibered $\proj^3$ of $\proj_F(N)$ over
$[F_1]$ was intrinsically described in section (\ref{resmap}).
We see
$$\proj^3= \proj(Hom_{\com}(K, Sym^3(V^*))) \eqq \proj(Sym^3(V^*))$$
where $V^*$ has a basis given by $(Z_1/Z_0)$ and $(Z_2/Z_0)$ with
induced $\com^*$-weights $-w_1+w_0$ and $-w_2+w_0$ respectively.
For weights $\barr{w}$ such that no two sum to twice
the third, the $\com^*$ action on the
fibered $\proj^3$ has $4$ isolated fixed points. The total
number of fixed points for the $\com^*$-action on $X(3)$ is
$31=1+12+6+3\cdot 4$.

        The $\alpha$-weights at the $4$ fixed points
in the fibered $\proj^3$ are obtained in the following manner.
First the differential of the $\com^*$-action on $H(3)$ at $F_1$
is determined. Then the blow-up along the locus $F$ is examined.
The sequence (\ref{xact}) of section (\ref{resmap})
contains all the necessary data.
The four fixed points in $\proj(Sym^3(V^*))$ are
$P_{r,s}=[(Z_1/Z_0)^r(Z_2/Z_0)^s]$ for non-negative integers $r$, $s$
with sum  $r+s=3$.

We tabulate the weight formulas for the four points $P_{r,s}$.
First the six  $\alpha$-weights:
\begin{equation*}
\begin{array}{lllllll}
P_{3,0}: & w_1-w_0 & w_2-w_0 & 2w_2-w_1-w_0 & w_1-w_2 & 2w_1-2w_2
& 3w_1-3w_2 \\
P_{2,1}: & w_1-w_0 & w_2-w_0 & w_2-w_0 & w_2-w_1 & w_1-w_2 & 2w_1-2w_2 \\
P_{1,2}: & w_1-w_0 & w_2-w_0 & w_1-w_0 & 2w_2-2w_1 & w_2-w_1 & w_1-w_2 \\
P_{0,3}: & w_1-w_0 & w_2-w_0 & 2w_1-w_2-w_0 & 3w_2- 3w_1 & 2w_2-2w_1
& w_2-w_1
\end{array}
\end{equation*}
\noindent
The $\beta$-weights for the four point $P_{r,s}$ all include the
the six weights:
\begin{equation*}
\begin{array}{lll}
-dw_0 & & \\
-(d-1)w_0-w_1 & -(d-1)w_0-w_2 & \\
-(d-2)w_0-2w_1 & -(d-2)w_0-w_1-w_2 & -(d-2)w_0-2w_2
\end{array}
\end{equation*}
The additional three $\beta$-weights at the points $P_{r,s}$ are:
\begin{equation*}
\begin{array}{llll}
P_{3,0}: & -(d-3)w_0-2w_1-w_2 & -(d-3)w_0-w_1-2w_2 & -(d-3)w_0-3w_2 \\
P_{2,1}: & -(d-3)w_0-3w_1 & -(d-3)w_0-w_1-2w_2 & -(d-3)w_0-3w_2 \\
P_{1,2}: & -(d-3)w_0-3w_1 & -(d-3)w_0-2w_1-w_2 & -(d-3)w_0-3w_2 \\
P_{0,3}: & -(d-3)w_0-3w_1 & -(d-3)w_0-2w_1-w_2 & -(d-3)w_0-w_1-2w_2
\end{array}
\end{equation*}
As before, the $\alpha$ and $\beta$-weights at the
fixed points on the other fibered $\proj^3$'s can be obtained
by permuting the weights $\barr{w}$ in the above formulas.
The Bott Residue Formula now yields:
$$c_6(E(3,d))={9\over2}d^6-27d^5+{9\over 2}d^4+{423\over 2}d^3-229d^2
-{829\over 2}d+525.$$
The algebraic computation was done on MAPLE with weights $\barr{w}=
(0,1,3)$.

\noindent
Department of Math, Harvard University,
harris@@math.harvard.edu

\noindent
Department of Math, University of Chicago,
rahul@@math.uchicago.edu

\end{document}